\documentclass[preprint,12pt]{elsarticle}
\usepackage{amsmath, amssymb, amsthm}
\usepackage{hyperref}
\usepackage{enumitem}
\theoremstyle{plain}
\newtheorem{lemma}{Lemma}
\newtheorem{theorem}{Theorem}
\newtheorem{example}{Example}
\theoremstyle{remark}
\newtheorem*{remark}{Remark}

\begin{document}

\begin{frontmatter}

\title{A Countable, Dense, Dedekind-Complete Subset of $\mathbb{R}$ Constructed by Extending $\mathbb{Q}$ via Simultaneous Marking of Closed Intervals with Rational Endpoints}

\author{Slavica Mihaljevic Vlahovic$^a$}
\author{Branislav Dobrasin Vlahovic$^b$\corref{cor1}}
\ead{vlahovic@nccu.edu}
\address{$^a$University of Zagreb, 41000 Zagreb, Croatia}
\address{$^b$North Carolina Central University, Durham, NC 27707, USA}
\cortext[cor1]{Corresponding author}

\begin{abstract}
This article explores the model-dependent nature of set cardinality, emphasizing that it is not an absolute measure. The cardinality of a set may differ across models, even when those models are internally consistent and complete within their respective universes. This inherent relativity is a foundational concept in modern set theory and underlies why cardinality, and the Continuum Hypothesis (CH), are intrinsically tied to the framework of the chosen model. For example, while the diagonal argument shows that the real numbers are non-denumerable within the Zermelo-Fraenkel (ZF) framework, the precise cardinality of the continuum remains unsettled and varies across different extensions of ZF. Under the inner-model axiom \textit{V = Ultimate L}, CH holds; in contrast, it is refuted under Martin’s Axiom. There is also a multiverse view, in which the continuum hypothesis is true in some models and false in others, all of which are considered mathematically valid. The Löwenheim–Skolem theorem reinforces this view by showing that no first-order theory can have only non-denumerable models. Thus, notions such as “denumerable” and “non-denumerable” become relative within any axiomatic framework, necessitating the relativization of cardinal numbers, as originally suggested by Skolem and Carnap.

To examine these issues, we construct and analyze two analogous sets: one within the ZFC (ZF plus the Axiom of Choice) framework, and another within Wang’s $\Sigma$ model, which differs from ZFC. Both constructions follow Gödel’s idea of building larger sets from smaller ones.

The first set, \(S_m\) is constructed from a countable collection \(S_Q\) of all closed intervals with rational endpoints, arranged in a sequence $\{S_Q\}$. All intervals are pre-marked through a simultaneous marking process, which assigned irrational points to all intervals in one unified step, respecting the structure of the sequence $\{S_Q\}$ and all its subsequences. Each interval is used exactly once, and no two distinct irrational numbers are assigned to the same interval.  While each interval \([a_n, b_n]\) contains infinitely many subintervals \([a'_n, b'_n] \subset [a_n, b_n]\), each of which is assigned a marked number, the construction ensures that exactly one marked number is assigned to each interval \([a_n, b_n] \in S_Q \). 
All intervals are marked without omission, ensuring that in every interval from the begining, as the outset, is one marked or assigned irrational number. Such marking reults in the set \(S_m\), which contains all rational endpoints from \(S_Q\) and the assigned irrational numbers. 

Within ZF, $S_m$ exhibits two critical properties: it is everywhere dense in \( \mathbb{R} \), and as it will be shown it is Dedekind complete. Since \(S_Q\) contains all closed intervals with rational endpoints, it enables inside \(S_Q\) identification of any Cauchy sequence representing any real number. 
The Axiom of Choice plays a critical role here, as it ensures that there exists a global selection function that simultaneously assigns marked numbers to all intervals in \( S_Q \) respecting the strucute of the sequence $\{S_Q\}$ and its subsequences, and for the identification of intervals forming any Cauchy sequence, without explicitly constructing the sequence itself. Because all intervals in such identified sequences are pre-marked from the outset, and that marking respects the nested structure of all intervals, that in all intervals there are numbers simultaneously, ensures that their intersection remains nonempty and includes the irrational number represented by the sequence. This leads to the remarkable conclusion that \( S_m \) is both complete and countable: a countable set of irrational numbers is distributed across a countable set of intervals, one per interval.

The second set is constructed within Wang’s $\Sigma$ model by systematically inserting irrational numbers between rational numbers using a diagonal procedure. Unlike in \( S_m \), these irrational numbers are not arbitrarily assigned, but are instead explicitly defined through successive diagonal constructions. The process begins with all rational numbers forming the \( 0 \)-th layer. The first layer is formed by applying a diagonalization procedure to this layer, thereby generating new irrational numbers. Subsequent layers are built by recursively applying similar diagonal procedures to the union of all previous layers. The process concludes after \( \omega \) steps, when further subdivision is no longer possible, consistent with the atomic nature of real numbers. The resulting set is equivalent to \( S_m \), and supports the conclusion that constructively defined real numbers are countable, apparently challenging the traditional interpretation of the diagonal argument. However, it is important to note that this construction applies only to constructible reals and excludes non-constructible ones. Nevertheless, the contradiction remains significant, as the diagonal argument is typically applied to explicitly definable numbers (i.e., those given by digital sequences).

Further analysis in the paper identifies tensions between foundational principles in set theory, especially between the diagonal argument and the Nested Interval Property (NIP). These issues are not limited to one argument, but recur across multiple classical proofs of uncountability. While these results may resonate with constructivist perspectives, especially those critical of the Axiom of Choice, they also motivate a broader reevaluation of the countability of the real numbers within formal mathematical frameworks.
\end{abstract}


\end{frontmatter}

\section{Introduction}

G\"{o}del conceived a small and constructible model universe denoted \emph{L}), generated by iteratively building larger sets from the empty set. In 1938, he put forward a new axiom for set theory: the axiom \emph{V=L}, asserting that every set is constructible.
By considering the inner model in which \( \mathsf{V} = \mathsf{L} \) holds, G\"{o}del demonstrated the relative consistency of the Zermelo–Fraenkel axioms (ZF), ZF augmented with the Axiom of Choice (ZFC), and the Generalized Continuum Hypothesis (CH) \cite{Godel}. While he initially entertained the idea of adopting \( \mathsf{V} = \mathsf{L} \) as a foundational axiom, he later rejected it, believing the axiom to be false despite its consistency within ZF. In his later view, the continuum hypothesis itself was not true, and not every uncountable set of real numbers has the same cardinality as \( \mathbb{R} \). He conjectured instead that \( 2^{\aleph_0} = \aleph_2 \) \cite{Godel1,Godel2}. As he wrote, \emph{``One may on good reason suspect that the role of the continuum problem in set theory will be this: that it will finally lead to the discovery of new axioms which will make it possible to disprove Cantor's conjecture''} \cite{Godel3}.

The independence of CH from ZFC was established through G\"{o}del's 1940 proof that CH cannot be disproven from the ZFC axioms, followed by Cohen’s 1963 result that it also cannot be proven within them \cite{DanaScott}. This led to ongoing efforts to introduce new axioms to ZFC that might resolve the continuum problem. Two of the most prominent candidates are the Ultimate-\( \mathsf{L} \) axiom, which is an all-encompassing generalization of G\"{o}del’s constructible universe consistent with CH and the cardinality \( \aleph_1 \); and Martin's Maximum axiom\cite{Martin,Solovay}, which asserts that anything conceived of using any forcing construction is a true mathematical entity, so long as the construction satisfies a certain consistency condition, compellingly pointing to $\aleph_2$ cardinality. 

Forcing axioms provide a method for extending the universe \( V \) of all sets through an iterative process of generic extensions: 
%\vskip-0.2cm
	\begin{equation}
V \subset V[g_0] \subset V[g_1] \subset \cdots \subset V[g_{\alpha}],
\label{Forcing}
\end{equation}

where each \( V[g_\alpha] \) represents a new universe generated by forcing over the previous one. Cohen’s method of forcing allows the creation of arbitrarily many new reals, beyond those of the model, resulting in extensions where the continuum has cardinality \( \aleph_2, \aleph_{315} \), or other transfinite values. However, this procedure produces a multiplicity of universes, a cloud of univrses, each with a different structure of the continuum, without a canonical way to distinguish which is virtual and which is real,
or even if any reflects the actual existence of mathematical objects.

To address this ambiguity, additional forcing axioms such as Martin’s Maximum and Woodin’s \( (*) \) axiom have been introduced \cite{Magidor}. These axioms enable the formulation of statements like \emph{``for all \( X \), there exists \( Y \) such that \( Z \)''} when referring to properties within a universe of sets, thus promoting certain existential principles that guide set-theoretic constructions.

However, there are inconsistencies between Martin's and (*) axioms and between these axioms and ZFC \cite{Cohen,Moore,Woodin}. For example, the statement entailed by \( (*) \), \emph{``for all sets of \( \aleph_1 \) reals, there exist reals not in those sets''}, clearly contradicts CH. Recent work by Aspero and Schindler \cite{Aspero} shows that a strong version, Martin’s Maximum\(^{++} \), implies Woodin’s \( (*) \), further linking two of the most influential axioms. Both imply that the continuum has cardinality \( \aleph_2 \), yet both also highlight unresolved tensions within the set-theoretic landscape.

Woodin has proposed further axioms \( (*)^+ \) and \( (*)^{++} \), which quantify over the entire power set of the reals, extending the expressive strength of \( (*) \). However, these stronger axioms are known to conflict with Martin’s Maximum in certain models \cite{Woodin}, questioning the existential statements about sets of reals in Martin's maximum axiom, and if if the reals can be coherently maintained in a unified set-theoretic framework.

The problem can be further underlined by L\"{o}wenheim-Skolem theorem or the intimately related G\"{o}del completeness theorem and its model-theoretic generalizations.  These results imply that any first-order theory formulated in a countable language must have a countable model, that a set can be {\textquoteleft nondenumerable\textquoteright} in this relative sense and yet be denumerable {\textquoteleft in reality\textquoteright}. As Skolem sums it up, \emph{{\textquotedblleft even the notions {\textquoteleft finite\textquoteright}, {\textquoteleft infinite\textquoteright}, {\textquoteleft simply infinite sequence\textquoteright} and so forth turn out to be merely relative within axiomatic set theory\textquotedblright}} \cite{Bays}.  Thus, the notion of countability is not absolute, but model-relative.

One of the consequences of the theorem is the existence of a countable model of set theory, which nevertheless must satisfy the sentence asserting that the real numbers are uncountable. The diagonal argument requires that $\mathbb{R}$ is nondenumerable in all models, that ZF set theory has only nondenumerable models. However, by the L\"{o}wenheim-Skolem theorem this is impossible, no theory can have only nondenumerable models. If a theory has a nondenumerable model, it must have denumerably infinite ones as well.   
Any model of ZFC that contains uncountable sets must admit a countable model (from the perspective of a larger meta-theory) that satisfies the same first-order properties. 
%One striking consequence of the Löwenheim–Skolem theorem is that any model of ZFC that contains uncountable sets must admit a countable model (from the perspective of a larger meta-theory) that satisfies the same first-order properties. Hence, although the diagonal argument asserts that \( \mathbb{R} \) is uncountable, ZF cannot preclude the existence of countable models in which the reals appear uncountable. In other words, if a theory admits an uncountable model, it must also admit denumerable models.

%This insight complicates the usual interpretation of the diagonal argument. The argument employs predicates definable only in a countable language, generating at most \( \aleph_0 \) subsets of the integers—thus enabling only countable characterizations of uncountability. The resulting enumeration of the continuum arises from an external meta-theoretic vantage point, not within the formal system itself \cite{Carnap}. While one can define bijections between \( \mathbb{R} \) and \( \mathbb{N} \) from outside a given model, such mappings cannot be constructed within the system if the system interprets \( \mathbb{R} \) as uncountable.

The non-denumerability is in contrast with predicates used in the diagonal argument, available only in a denumerably infinite multitude and yield no more than $\aleph_0$ subsets of the set of all integers, i.e., an enumeration of the continuum from outside the axiomatic system \cite{Carnap}.  The enumeration from outside utilizes the structure of the system as a whole, which cannot be reached by operations within the system. There are one-to-one correspondences between $\mathbb{R}$ and $\mathbb{N}$, but they all lie outside the given model. 

Consequently, the uncountability of the reals becomes a model-relative property. A set that is countable in one model may be uncountable in another. For every non-denumerable model \( M \) of set theory, there exists a countable model \( M' \subset M \) satisfying the same first-order axioms. Variables in \( M' \) range over a proper subset of those in \( M \), allowing a denumerable interpretation of otherwise uncountable sets. Therefore, it is not possible to rule out denumerable interpretation. Infinite sets may have different cardinal numbers and be nevertheless syntactically of the same cardinal number \cite{Carnap}. In a given axiomatic set theory $M$, there always exists a more comprehensive theory in which all infinite sets of $M$ prove equivalent, namely denumerable. Every non-denumerable set becomes denumerable in a higher system or in an absolute sense.

%These results provide a foundational basis for reexamining the cardinality of the continuum. Rather than taking the uncountability of the real numbers as an absolute conclusion, one may consider it a relative feature, contingent on the model and the axioms adopted. The remainder of this paper presents two explicit constructions of countable sets with properties commonly attributed to \( \mathbb{R} \), offering insight into the model-dependence of cardinality and the interpretational flexibility inherent in set theory.

\section{Discussion}

Let us again consider  G\"{o}del's axiom $V = L$ that all sets are constructible and examine the impact of the L\"{o}wenheim-Skolem theorem. As shown in \cite{Barwise} the theorem requires that ZF plus $V = L$ has an $\omega$-model which contains any given set of real numbers. It follows from the statement that for every real $s$, there is a countable $M'$ such that $M'$ is an $\omega$-model for ZF plus $V = L$ and $s$ is represented in $M'$, which follows from the theorem. So, by the L\"{o}wenheim-Skolem theorem, a model containing $s$ can satisfy {\textquoteleft $s$ is constructible\textquoteright} (because it satisfies {\textquoteleft $V = L$\textquoteright}, and {\textquoteleft $V = L$\textquoteright} says everything is constructible) and be an $\omega$-model. In the examples below, we will construct such $\omega$-models. 

Let us first follow the L\"{o}wenheim-Skolem theorem and expand ZF model with AC. 
While AC is independent of the ZF set theory, which means that both the axiom itself and its negation are consistent with ZF, AC expands ZF. 
It allows the construction of set $S_m$, which is 
%and to recognize in the set $S_m$ every possible Cauchy sequence representing any number, even if no sequence is ever explicitly defined,
%It also allows the construction of any sequence of any number from the intervals with rational endpoints that are present in the set, similarly to 
some kind of a variation of the {\textquoteleft $V = L$\textquoteright} approach.  
As described in Theorem 1, set $S_m$ is generated by marking or identifying an arbitrarily selected irrational number within each interval with rational endpoints. It will be shown that such a countable set is everywhere dense in $\mathbb{R}$, and Dedekind complete,  i.e., that it has properties of set $\mathbb{R}$. 

\subsection{\bf Main Construction and Theorem}
Let us start with a set of rational numbers and, in a constructive way under the ZFC formalism, create a set equivalent to the set of real numbers.

\begin{theorem}[Construction and properties of the set $S_m$]
Enumerate all closed intervals with rational endpoints:
\begin{equation}
S_Q = \{ I_n = [a_n, b_n] : a_n, b_n \in \mathbb{Q},\; a_n \le b_n,\; n\in\mathbb{N} \}.
\label{SQ}
\end{equation}
Suppose each interval $I_n$ is assigned exactly one irrational mark $x_n \in I_n$ under the following structural marking axioms:
\begin{enumerate}[label=(S\arabic*),leftmargin=3em]
\item \textbf{Coverage:} Every interval $I_n$ receives exactly one irrational mark:
\begin{equation}
\forall n \in \mathbb{N},\quad x_n \in I_n.
\label{xninIn}
\end{equation}
\item \textbf{Nested Interval Marking:} Every subinterval with rational endpoints of a marked interval,  \([a',b']\subseteq I_n\),  is also marked from the outset. Thus, no interval or subinterval, regardless of its depth or nesting, remains unmarked. This implies that each \(I_n\) contains infinitely many marked points.

\item \textbf{Global Simultaneous Marking:} All intervals and subintervals are marked simultaneously in a global, structurally coherent manner. Specifically, the marking does not depend on sequential adaptation but is preassigned once and for all, fully respecting the nested and overlapping structure of all its subintervals.
\item \textbf{Universal Markability:} For every interval $I_n=[a_n,b_n]$,  any point $x_n\in I_n$ is eligible to be assigned as a mark $x_n\in M$.
\end{enumerate}
Let the set of all such marks be:
\begin{equation}
M = \{x_n : n \in \mathbb{N}\} \subset \mathbb{R} \setminus \mathbb{Q}.
\label{SetM}
\end{equation}
and structurally marked set
\begin{equation}
S_m := \mathbb{Q} \cup \{x_n : n \in \mathbb{N}\} = \mathbb{Q} \cup M.
\label{SetSm}
\end{equation}
Then the set $S_m$ is countable, everywhere dense in $\mathbb{R}$, and Dedekind complete, that is, every Dedekind cut in $S_m$ is realized by an element of $S_m$.
\end{theorem}

\vskip 0.3cm
{\it Proof:}
\vskip 0.3cm

\subsection*{\bf Countability}
$\mathbb{Q}$ and $M$ are both countable, so $S_m$ is a countable union of countable sets.

\subsection*{\bf Density}
Since $\mathbb{Q} \subset S_m$ and $\mathbb{Q}$ is dense in $\mathbb{R}$, it follows that $S_m$ is dense in $\mathbb{R}$.

\subsection*{\bf Dedekind Completeness}
To prove Dedekind's completeness, let us first consider this lemma:
\begin{lemma}[Dedekind Completeness of $S_m$ via Gap‐Point \& Squeeze]
Under the structural marking axioms (S1)--(S4) of Theorem 1, the set $S_m = \mathbb{Q} \cup M$ is Dedekind complete: every Dedekind cut in $S_m$ is realized by an element of $S_m$.
\end{lemma}

\begin{proof}
Let $A \mid B$ be a Dedekind cut in $S_m$, and set
\begin{equation}
t = \sup(A \cap \mathbb{Q}) = \inf(B \cap \mathbb{Q}) \in \mathbb{R}.
\label{Cutt}
\end{equation}
We show $t \in S_m$.

\medskip\noindent\textit{Case 1: $t \in \mathbb{Q}$.} 
Then $t \in \mathbb{Q} \subset S_m$, so the cut is filled.

\medskip\noindent\textit{Case 2: $t \notin \mathbb{Q}$.} 
We must prove $t \in M$.

\smallskip\noindent\emph{(a) Nested‐interval representation.} 
Since $t$ is the gap‐point in $\mathbb{Q}$, choose
\begin{equation}
a_1 < a_2 < \cdots < t < \cdots < b_2 < b_1, \quad b_n - a_n < \tfrac1n,
\label{Sequencesanbn}
\end{equation}
and let $I_n = [a_n,b_n]$. Then
\begin{equation}
I_1 \supseteq I_2 \supseteq \cdots, \quad \text{and} \quad \bigcap_{n=1}^\infty I_n = \{t\}.
\label{Insubset}
\end{equation}

\smallskip\noindent\emph{(b) Define $\alpha,\beta$.} 
Let
\begin{equation}
\alpha = \sup\{x \in M : x < t\}, \qquad \beta = \inf\{x \in M : x > t\}.
\label{alphabeta}
\end{equation}
Suppose that $t \notin M$, we only know $\alpha \le t \le \beta$.

\medskip\noindent\emph{Subcase 2a (strict gap $\alpha < \beta$).} 
Pick rationals $r^- < t < r^+$ with $\alpha < r^- < t < r^+ < \beta$. For large $n$,
\begin{equation}
I_n \subset (r^-, r^+),
\label{r}
\end{equation}
and by definition no mark lies in $(r^-, r^+)$, so $I_n$ would be unmarked—contradicting (S1). Thus $\alpha < \beta$ is impossible.

\medskip\noindent\emph{Subcase 2b (cluster $\alpha = \beta = t$).} 
Here every mark in $M$ accumulates at $t$. Now:

1. By (S1)--(S3), each $I_n$ is marked, so pick $x_n \in M \cap I_n$. Thus, for each $n$, both $x_n \in M$ and $t \in I_n$.

2. Since $\mathrm{diam}(I_n) \to 0$ and $t \in I_n$, we have $|x_n - t| \le \mathrm{diam}(I_n) \to 0$, hence $x_n \to t$.

3. For any fixed $x_k \ne t$, set $\delta_k = \tfrac12 |x_k - t| > 0$. Then for all $n$ with $\mathrm{diam}(I_n) < \delta_k$, we have $I_n \subset (t - \delta_k, t + \delta_k)$, so $x_k \notin I_n$. Hence no $x_k \ne t$ can lie in every $I_n$. The Nested Interval Property ensures that all $x_k \ne t$ are eventually excluded, and only $t$ survives.

4. Yet each $I_n$ must contain \emph{some} mark. Axiom (S2) requires that there are marked points $x \in M$ in every subinterval that, together with $t$, remain in every subinterval, regardless of how small, all the way to the intersection. If no mark survives through all intervals together with $t$, some deep subintervals would be left unmarked, contradicting (S2).

Thus, the only point able to survive \emph{all} levels is $t$ itself. Therefore, one of the $x_n$ is forced by the structure of the intervals not only to converge to $t$ but to be equal to $t$, which is allowed by (S4), and hence \(t\in M\).

In both subcases, we reach the conclusion that $t \in M \subset S_m$, contradicting the assumption that $t \notin S_m$. Thus, every Dedekind cut in $S_m$ is filled, and $S_m$ is Dedekind complete.
\end{proof}

%So, for any Dedekind cut in \(S_m\) with cut point \(t\in \mathbb{R}\), we have:

%\paragraph{Case 1: \(t \in \mathbb{Q}\).} Then \(t \in S_m\), and the cut is realized.
%{\it Case 1}: \(t \in \mathbb{Q}\). Then \(t \in S_m\), and the cut is realized.

%\paragraph{Case 2: By Lemma~\ref{lemma:limit}, \(t \in M\), so the cut is filled.
%{\it Case 2}: By Lemma 1. \(t \in M\), so the cut is filled.

%By Lemma~\ref{lemma:limit}, \(t \in M\), so the cut is filled.

%Therefore, there is no gaps in no gaps in the set $S_m$ and teh set $S_m$ is complete.
%%%%%%\end{proof}

\textbf{Conclusion of Theorem 1} 

By combining the results, it follows that \( S_m \) satisfies all the stated properties. 
Specifically, \( S_m \) is everywhere dense in \( \mathbb{R} \), Dedekind complete, and countable.  Therefore, \( S_m \) represents a rigorous construction of a set that is countable yet exhibits the properties of the real numbers. This proves the main Theorem 1.  
\hfill $\square$

This outcome is unexpected, but not comppletely, since for any arbitrarily selected numbr $r$ we can find in the set $S_Q$ collection of intervals $[a_n,b_n]$, every contining marked numbers, that are closing on $r$ and after an infinite overlapping will include only one marked numer that must be $r$.  

To demonstrate that select a number $r$ and recognize that the set \( S_Q \) consisting of intervals with rational endpoints is dense in \( \mathbb{R} \), for any \( \epsilon > 0 \), there exist intervals \([a_n, b_n] \in S_Q\) that can be 
iteratively selected  with smaller \( \epsilon_n \to 0 \) as \( n \to \infty \), constructing a sequence of nested intervals \([a_n, b_n]\) satisfying:
	\vskip-0.2cm
	\begin{equation}
	r - \epsilon \leqslant a_n \leqslant r \leqslant b_n \leqslant r + \epsilon.
	\label{Minterval}
	\end{equation}
	Therefore in the set \(S_Q\) it is always possible to recognize and select an infinite sequence of nested intervals  \(S_{Q_n}\), or better say in set \(S_Q\) there exists an implicit, infinite sequence of nested nonempty intervals  \(S_{Q_n}\) with endpoints that  close on \(r\) and uniquely select the point \(r\), as:       
	\begin{equation}
	a_n\leqslant r\leqslant b_n, \mbox{ for all } n. 
	\label{Con}
	\end{equation}

Dedekind considered a similar problem. The  number \(r\) can be defined through Dedekind's cut \(A\mid B\), with \(A\) and \(B\) as the disjoint components of the cut. The endpoints of the nested intervals \([a_n, b_n]\) serve as elements defining the cut \(A\mid B\), with \(a_n\) in component \(A\) and \(b_n\) in component \(B\). As \(n\) increases, the gap between \(a_n\) from \(A\) and \(b_n\) from \(B\) diminishes, reducing the distance between pairs \(a_n\) and \(b_n\) toward zero. This progressive narrowing necessitates that the only number fitting between every pair \(a_n\) and \(b_n\) is \(r\), the number that represents the gap created by the cut, as it was shown by Dedekind \cite{Dedekind}.

Since Theorem 1. is an important result, let us derive it independently. 

%\end{proof}
 %xxxxxxxxxxxxxxxxxxxxxxxxxxxxxxxxxxxxxxxxxx   

\section{Alternative proof for Completeness of the set $S_m$}
%\begin{theorem}[Te set $S_m$ does not have any gaps]
%%%%%%%%%%%%%55ABCD
%\end{proof}

\begin{theorem}[Dedekind Completeness of \( S_m \) via Dedekind Cuts]
\vskip 0.3cm

The set \(S_m\) as defined in Theorem 1, 
 is Dedekind complete. Any Dedekind cut \(A\mid B\) of  \(S_m\) into two disjoint components \(A\) and \(B\) does not create a gap. Specifically, the first component  \( A \)  has a greatest element, or the second component \( B \)  has a least element, or both. This means the set \(S_m\) as defined in Theorem 1, contains all real numbers, there is no a number that is not marked. Any new different marking will not produce marked numbers that are not already included in the marked set $S_m$. 

\end{theorem}
	
\begin{proof} 

Let  \(S_Q\) and \( S_m \) be, as defined in theorem 1.  The set \( S_Q \) consists of all closed intervals \([a_\alpha, b_\alpha]\) with rational endpoints \( a_\alpha, b_\alpha \in \mathbb{Q} \), \( a_\alpha \leq b_\alpha \). The set \( S_m \) is constructed by marking one irrational number \( x_\alpha \) in each interval  \( [a_\alpha, b_\alpha] \in S_Q \), under the conditions specified in the theorem, and includes all such numbers as well as  all  rational numbers \( \mathbb{Q} \). 
\vskip 0.2cm	

\textit{Relationship Between Cuts in \(\mathbb{Q}\) and \(S_m\)}

Marking an irrational number \(r\in S_m\) in every interval \([a_\alpha,b_\alpha]\in S_Q\) corresponds  to  creating Dedekind cuts \(A\mid B\) in \(\mathbb{Q}\).  Specifically, each marked \( r \), divides the rational numbers \( \mathbb{Q} \) into two disjoint subsets:  \( A = \{q \in \mathbb{Q} : q < r\}, \quad B = \{q \in \mathbb{Q} : q \geq r\}. \)  This partition satisfies the properties of a Dedekind cut: \( A \cup B = \mathbb{Q} \), \( A \cap B = \emptyset \), and for all \( a \in A \), \( b \in B \), \( a < b \).  

So, marking an irrational number in an \(S_Q\) interval is equivalent to representing a gap in \(\mathbb{Q}\) within that interval. 
%Similarly, as it will be shown, every cut resulting in a gap in \(\mathbb{Q}\) within an \(S_Q\) interval, is represented by a corresponding marked irrational number in \(S_m\).

\vskip 0.2cm	
\textit{Demonstrating Completeness of \(S_m\)}

%To show that \(S_m\) is Dedekind complete, consider any Dedekind cut \( A/B \) of \( S_m \). The proof focuses on cases where \(A\) lacks a greatest element or \( B \) lacks a least element, as these correspond to potential gaps. Cuts that do not create gaps correspond to elements that are already in \(S_m\).

After all intervals $I_n$ are marked and set $S_m$ is composed, consider inspecting that all irrational numbers are marked by performing an additional, new arbitrary marking of already marked intervals $I_n$ and check if such additional markig will producen a different marked number $t$ that is already not icluded in the marked set $S_m$. Such additional marking will split all numbers in the interval, including marked numbers, into two disjoint components $A$ and $B$. So, it represents a Dedekind cut 
\[ 
A\mid B
\]
in an interval $I_n$ and also a cut in the set of marked numbers $S_m$ separating them also in the components $A$ and $B$.

Since any marked number by the theorem is irrational, the new marked number $t$ also makes a Dedekind cut in rational numbers $\mathbb{Q}$:
\begin{equation}
t = (A\cap \mathbb{Q})\mid (B \cap \mathbb{Q}).
\label{eq2}
\end{equation}
%where \(t\) separates rational numbers in \( \mathbb{Q}\). 

As a marked number $t$ belongs to the set of marked numbers $S_m$, and as a marked element of the interval $I_n$, it must belong to one of the components $A$ or $B$ of the cut in $I_n$, therefore it belongs to either of the components $A$ or $B$ of both $S_m$, and $I_n$. 
%It will be shown that if \(t\) is in \(A\) that it is then the greatest element of component \(A\), and if it is in \(B\) that it is the least element of component \(B\). This means that cuts in \(S_m\) do not generate gaps in set \(S_m\).

It will be shown that this new marking does not create a gap in the set $S_m$ of marked numbers. 
Specifically, the first component \( A \) of marked numbers has a greatest element, or the second component \( B \) has a least element, or both. So, it will be shown that any new marking generates a marked numbers that already exist in set $S_m$.

Let us show that if \( t \in A \), then it is the greatest element of component \( A \). Suppose \( t \) is not the greatest element of component \( A \). Then there exists a 
\begin{equation}
t' \in A \mbox{ such that } t<t' 
\label{eq4S}
\end{equation}
The \( t' \) cannot be a rational number because \( t < t' \) and, by (\ref{eq2}) the definition of the cut \(t\), it would require \(t' \in B \cap \mathbb{Q}\), contradicting \( t' \in A \).

If \( t' \) is an irrational number, it represents a cut \( C\mid D \) in \( \mathbb{Q} \), which creates a gap in \( \mathbb{Q} \), and for this reason, the first component \(C\) of that cut \(t'\) does not have the greatest element. Because \( t < t' \), we have 
\( (A \mid B) < C\mid D \), 
implying \( (A \cap S_m) \subset C \). If \( t'' \) is any marked number from \( C \setminus A \) (since \(C\) does not have the greatest element and \( S_m \) is everywhere dense), then \( t < t'' \) would require \( t'' \) to be in \( B \). 
Contrarily, because \(t'' \leqslant t'\) (as \(t'\) represents a cut in \(\mathbb{Q} \) and is in \(A\), as referenced in (\ref{eq4S})), it requires that \(t'' \in A\).

However, having \(t''\) in both \(A\) and \(B\) is not possible because \(A\) and \(B\) are disjoint components of the cut. 
This contradiction implies there is no \( t' \in A \) such that \( t < t' \), and so \( t \) must be the greatest element of \( A \).

Similarly, if \( t \in B \), it can be shown that \( t \) is the least element of component \( B \). 

This ensures that \( S_m \) satisfies Dedekind completeness: any cut in \( S_m \) does not create a gap or require additional elements.
Since \( t \) is either the greatest element of \( A \) or the least element of \( B \), it follows that any Dedekind cut \( A\mid B \) in \( S_m \) does not create a gap. Every cut is represented by a unique marked point \(t \in S_m\), ensuring that \( S_m \) satisfies Dedekind completeness. 

Thus, \( S_m \) is complete in the Dedekind sense, having no gaps, affirming the Theorem. 
\end{proof}

Let us show Dedekind complentess also on this independet way:

\begin{theorem}[Dedekind Completeness of \( S_m \) via Missing Elements]
\label{thm:completeness-missing}
The set \( S_m \), as defined in Theorem~1, is Dedekind complete. That is, \( S_m \) contains every real number: no real number is missing from \( S_m \).
\end{theorem}

\begin{proof}
Suppose, for contradiction, that there exists a real number \( t_x \notin S_m \). Then \( t_x \) determines a Dedekind cut in \( \mathbb{Q} \), and can be represented as the unique limit of a nested sequence of closed intervals with rational endpoints:
\[
I_n = [a_n, b_n], \quad a_n, b_n \in \mathbb{Q}, \quad I_1 \supseteq I_2 \supseteq \cdots, \quad \bigcap_{n=1}^\infty I_n = \{t_x\}.
\]

Since \( S_Q \) enumerates all intervals with rational endpoints, each \( I_n \in S_Q \), and the sequence \( \{I_n\} \) is a subsequence of this enumeration. By the global, simultaneous marking condition (S3), every interval and all its subintervals are assigned a mark from the outset.

Thus, for each \( n \), the interval \( I_n \) contains not just one marked point but infinitely many, inherited through the structural marking:
\[
M \cap I_n \neq \emptyset \quad \text{for all } n \in \mathbb{N}.
\]
Since the sequence \( \{I_n\} \) is nested and every \( I_n \) intersects \( M \), the finite intersections remain nonempty:
\[
\bigcap_{n=1}^m (M \cap I_n) \neq \emptyset \quad \text{for all } m \in \mathbb{N}.
\]
Now consider the infinite intersection. Because the intervals satisfy \( \bigcap_{n=1}^\infty I_n = \{t_x\} \), we obtain:
\[
\bigcap_{n=1}^\infty (M \cap I_n) = M \cap \left( \bigcap_{n=1}^\infty I_n \right) = M \cap \{t_x\}.
\]

To conclude \( t_x \in M \), we reason as follows:

Every interval \( I_n \) contains both the point \( t_x \in I_n \) and at least one marked point \( x_n \in M \cap I_n \). By condition (S3), the marking is structural and simultaneous, so no “gaps” are left that would prevent progression through the entire nested sequence.

Since \( t_x \) survives through all \( I_n \), and some mark exists in each interval together with $t_x$, structural marking ensures that at least one of these marks also persists through all intervals, i.e., there exists some \( x_k \in M \) that lies in every \( I_n \). All other marks $x_n\neq t_x$ will be eliminated with smaller intervals.  Thus to avoid leaving a deep interval unmarked, the limit point must have been one of  marks, 

But the only point that lies in every \( I_n \) is \( t_x \). Therefore, \( x_k = t_x \), and so \( t_x \in M \subseteq S_m \), contradicting our assumption.

Thus, \( t_x \) is a marked point, and every real number lies in \( S_m \), proving Dedekind completeness.
\end{proof}

%abcdaaaaaaaaaaaaaaaaaaaaaaaaaaaaaaaaaaaaaabbbbbbbbbbbbbbbbbbbbbbbbbbbbbbbbbcccccccccccccccccccccc
%\subsection*{Density}

\section{Explicit sequence of everywhere dense set of real numbers}

The properties of the set \(S_m\) clearly demonstrate that it supports the Löwenheim-Skolem theorem rather than the diagonal argument. This result, along with the Löwenheim-Skolem theorem, shows that our understanding of the nondenumerability of real numbers consists in our knowing what it is for this to be proved and not in our grasp of a model. 

The set \(S_m\), while complete and countable, does not provide an explicit sequence of all real numbers. This means that for a specific given number, it is not possible to identify its particular \(n\)-th position in the sequence. This limitation arises because the intervals in the set \(S_m\) even countable are not enumerated and to each interval is assigned one irrational number that is not specified. 
To address this, considered will be a set equivalent to \(S_m\), which is generated by inserting irrational numbers between rational numbers. Contrary to the previous example, this method specifies all generated irrational numbers at each step, thus providing a sequence in which numbers have clearly defined positions. This approach merges the concepts of Gödel's Constructible Universe (V = L) with the earlier proposed Wang's \(\Sigma\) model \cite{Wang}, which differs from ZF models. However, instead of using Wang's original work \cite{Wang}, let us modify the model introduced in \cite{Frankel}.

The model begins with a \(0\)-th layer, which consists of some denumerable totality of objects, which may be taken to be, for instance, the positive integers or all finite sets  built up from the empty set. In our next theorem, we will consider that the \(0\)-th layer comprises rational numbers. Then, the first layer includes all the objects of the \(0\)-th layer and, in addition, all those sets of these objects that correspond to conditions that contain bound variables ranging over objects of the \(0\)-th layer, but no bound variables ranging over the first or higher layers. Generally, the \(n+1\)-st layer contains all the objects of the \(n\)-th layer, along with all such sets of these objects determined by conditions whose bound variables range over objects of the \(n\)-th layer at most. The hierarchy of layers is continued beyond the finite ordinals. Layer $\omega$, for instance, is the sum of all finite layers, and layer $\omega + 1$ contains, in addition, also such sets as are determined by conditions whose bound variables range over the objects of layer $\omega$ at most. Since all the sets of layer $\alpha$ are enumerable by a function $E_{\alpha}$ ranging over entities of a layer, all sets of any subsystem $\Sigma_{\alpha}$ of $\Sigma$ are therefore enumerable. Therefore, this model requires that any infinite set can be enumerated in an appropriate partial system of $\Sigma$.  
We consider here Wang's model, but the conception of cumulative layers (or orders) is not original with Wang \cite{Godel4,Quine}, and there are also other attempts to realize it \cite{LAbbe}.   
	
In the next theorem, we present our model. Its \(0\)-th layer consists of a set of rational numbers. Layer one includes all numbers from the \(0\)-th layer and all irrational numbers created by simultaneously applying an infinite number of diagonal constructions on the \(0\)-th layer. The layer \(n+1\) is obtained by applying \(\omega\) diagonal constructions on layer \(n\), which also incorporates all previous layers. More details are provided below.

\begin{theorem}[Explicit sequence of everywhere dense and complete set of real numbers]

The construction described below generates an explicit sequence that contains all real numbers.
\end{theorem}
 Let us consider an interval of numbers \([M, W]\), starting by writing a denumerable sequence of all rational reals \(A_\nu\) that belong to that interval, as specified in the corresponding array (\ref{ar1}). Each number will be represented by an infinite binary sequence, even if it has a finite number of digits, with each digit \(a_{\nu\mu}\) being either 0 or 1:
\begin{equation}
{\begin{array}{clcr}
	A_1=(a_{11},a_{12},...,a_{1\mu},...)\\
	A_2=(a_{21},a_{22},...,a_{2\mu},...)\\
	\vdots
	\;\;\;\;\;\;\;\;\;\;\;\;\;\;\;\;\;\;\;\;\;\;\;\;\;\;\;\;\;\;\;\;\;\;\; \\
	A_{\nu}=(a_{\nu1},a_{\nu2},...,a_{\nu\mu},...)\\
	\vdots
	\;\;\;\;\;\;\;\;\;\;\;\;\;\;\;\;\;\;\;\;\;\;\;\;\;\;\;\;\;\;\;\;\;\;\;
	\end{array}}
\label{ar1}
\end{equation}
The sequence \(S_r\) contains only all rational numbers, and each rational number has a well-defined place within that sequence. Let us now use Cantor's diagonal construction to generate numbers not contained in this sequence. Such numbers are defined by sequences of digits \(d_{\nu1}, d_{\nu2}, \ldots, d_{\nu\mu}, \ldots\), where each \(d_{\nu\mu}\) is either 0 or 1, determined such that \(d_{\mu\mu} \neq a_{\nu\nu}\). From this sequence of \(d_{\nu\mu}\), we formulate the element \(D_1 = (d_{11}, d_{12}, \ldots, d_{\nu\mu}, \ldots)\), ensuring that \(D_1 \neq A_\nu\) for any index \(\nu\). Thus, a number different from all numbers in the sequence containing all rational numbers is created. Therefore, it must be an irrational number. This will be the first irrational number generated and will be incorporated into the sequence \(S_r\) as the first element, which shifts all elements in the original sequence from the \(n\)-th to the \(n+1\)-th place. This recursive construction will be repeated; it will be applied again on the newly generated set \(S_r\), generating in the same way the second irrational number and placing it in the sequence \(S_r\) in the first position, shifting all elements in \(S_r\) again, then the third irrational number, and so on.
All numbers generated in this way are well defined, all digits \(d_{\nu\mu}\) of such generated irrational numbers \(D_1, D_2, \ldots, D_\nu, \ldots\) are defined and known. This is because all digits in the initial sequence \(S_r\) that contains all rational reals are defined. In creating irrational numbers by Cantor's diagonal construction, there is no freedom because each diagonal element 0 must be switched to 1, or 1 to 0; one must apply \(d_{\nu\nu} \neq a_{\nu\nu}\) for each diagonal element. Since in the initial sequence \(S_r\) of rational numbers all \(a_{\nu\mu}\) are known and fixed, all numbers \(D_1\), \(D_2\), \ldots that are produced with the diagonal construction are defined.

To generate the first layer in our model, the construction will be modified so that in the first step, not only one but infinitely many irrational numbers are generated. This is possible because all decimal places \(a_{11}, a_{12}, \ldots, a_{1\mu}, \ldots\) can be used simultaneously and an infinite number of diagonal constructions can be applied simultaneously in the first step. There will be a diagonal construction that starts with the first decimal place \(a_{11}\), one that starts with the second decimal place \(a_{12}\), one with the third \(a_{13}\), and so on. New irrational numbers will be created using an infinite number of diagonal constructions, each starting from a different decimal place. After step one, the first layer, denoted as \(\Sigma_1\), will consist of an infinite sequence of \(\aleph_0\) rational numbers from the \(0\)-th layer and an infinite sequence of \(\aleph_0\) irrational numbers produced by an infinite number of diagonal constructions applied to the set of rational reals of the \(0\)-th layer. These two sets, the \(0\)-th layer and the irrational numbers generated by diagonal constructions, will be combined into one sequence by alternating the first number from the \(0\)-th layer with the first number from the list of created irrational numbers, and so on, to create the first layer, the \(\Sigma_1\) set. The construction from step one will be repeated, applying an infinite number of diagonal constructions to the sequence that represents the just created set \(\Sigma_1\), creating an additional \(\aleph_0\) irrational numbers that will then be combined with set \(\Sigma_1\) to create set \(\Sigma_2\). This recursive construction will continue indefinitely, creating an additional \(\aleph_0\) of newly created irrational numbers each time.

However, this construction will conclude after \(\omega\) steps due to the atomistic structure of real numbers, which prevents endless division of an interval. The process concludes exactly after \(\omega\) divisions because an infinite division of intervals results in singletons that cannot be divided further. At each step of creating layers, infinite amounts of irrational numbers are inserted between the algebraic numbers of the \(0\)-th layer, dividing each segment defined by a pair of rational numbers further until it becomes a singleton, as follows from the BW theorem or the NIP.
%%%%%%%%%%%%\end{theorem}

\begin{proof}

This Theorem is equivalent to Theorem 1. In both cases, irrational numbers are generated and inserted between rational numbers, subdividing the intervals into smaller and smaller segments. Decimal places \(a_{\nu\mu}\) that represent numbers in (\ref{ar1}) are equivalent to nested intervals in \(S_Q\), and making chnages in decimal places in  (\ref{ar1}) is equivalent to marking numbers in intervals \(S_Q\). In both cases generated is a set of everywhere dense irrational numbers. While in Theorem 1 completnes of set \(S_m\) arrises from its density and the fact that any number must be marked, and that there is no exclusion of numbers, here completeness follows from the fact that any missing number generated by diagonal argument is actullay part of the set of irrational numbers generated by this construction.  There is also a proof that any irrational number can be generated by a diagonal sequence  \(D_1, D_2, \ldots, D_\nu, \ldots\) in the set of rational numbers, by Gray, theorem 3 in \cite{25}. He proved that by making permutations in the positions of the elements of the initial rational sequence it is possible to obtain the condition that will generate the specific targeted irrational number. By our construction, we are not creating direct permutations of the positions of the elements in the rational sequence, but instead doing that indirectly by making in each step an infinite number of changes of the elements in the initial sequence, by adding to the initial sequence in each step an infinite number of new irrational numbers that have decimal places different from the initial sequence. Doing that using an infinite number of diagonals, changing all digits on each of the diagonals, we are actually in each step creating an infinite number of permuted initial sets considered by Gray. The permutations we are creating are predetermined by the initial set, but the fact that generated irrational numbers have a random distribution of the digits, ensures that conditions for generating any number can be achieved, i.e., that method does not exclude any number that one can specify. Any number that one can specify can be in principle generated, because the construction is going indefinitely, it is only a question in which step it will be produced. Actually, if one specifies all decimals of a number, it can be calculated at which sep that number will be generated.  

The diagonal argument cannot be used to create missing numbers not generated by this model because it is utilized to generate new numbers. Any \textquoteleft missing\textquoteright  number produced by diagonal construction is, in fact, just a number that belongs to this model. 
\end{proof}

To our knowledge, this is the first presentation of an explicit Wang's model.

\section{Remarks to the non-denumerability proofs}

There are many variations of proofs of the uncountability of real numbers. However, fundamentally, all of them share a common basis: they revolve around infinite nested sequences of intervals. We will discuss the critical point that appears in all these proofs, the interpretation of the intersection of these intervals and the nature of the number that resides within the intersection. According to the NIP, the intersection of an infinite nested sequence of closed nonempty intervals must be a zero-size interval, a singleton, or a single number. This number must be one of the elements from the intervals of the sequence, which, in all uncountability proofs, are the numbers of the supposed denumerable sequence assumed to include all numbers. Contrarily, in all nondenumerability proofs, that number is declared as the number that is outside the denumerable sequence, and its existence is used as the argument that the sequence does not contain all numbers. This discrepancy between the NIP, which posits that the number is a member of the sequence, and the uncountability proofs, which assert the opposite, will be elaborated upon, considering the proofs from 1873 and 1879, and the diagonal argument.

The proofs from 1873 and 1879 employ contradiction, and both begin by assuming that the real numbers in the interval \([\alpha, \beta]\) can be enumerated in a sequence:
\begin{equation}
                               \omega_1, \omega_2, ... , \omega_\nu, ...
\label{eq13}
\end{equation}  
and then establish a contradiction by demonstrating that there exists a number \(\eta\) that is not included in (\ref{eq13}). The 1873 article states that \emph{`` \(\eta\) can be any number in the interval \([A, B]\), where \(A\) and \(B\) are limits \([\omega_{a\nu}, \omega_{b\nu}]\) when \(\nu=\infty\)",} which is impossible because that is a zero-size interval, essentially a single number, and therefore no number \(\eta\) can reside within that interval. The same assertion appears in the 1879 article. The conclusion is that \emph{`` the number \(\eta\) will be inside that interval, while all listed numbers will be outside"}, which is incorrect because no number can reside inside a single number. These are significant objections that warrant further elaboration.

Let us refute \textbf{the original proof from 1879} to illustrate the problem. Starting with the sequence (\ref{eq13}) and an arbitrary interval \([\alpha, \beta]\), where \(\alpha < \beta\), identified can be a real number \(\eta\) that does not occur in the sequence (as a member of it). 

Some of the numbers from the sequence (\ref{eq13}) definitely occur within the interval \([\alpha, \beta]\). Among these, let \(\omega_{k1}\) be the number with the smallest index, and \(\omega_{k2}\) be the number with the next larger index. Let the smaller of these two numbers, \(\omega_{k1}\) and \(\omega_{k2}\), be denoted by \(\omega_a^{1}\), and the larger by \(\omega_b^{1}\). These two numbers define the first closed interval \([\omega_a^{1}, \omega_b^{1}]\). Within this interval, there are no numbers \(\omega_\nu\) from sequence (\ref{eq13}) for which \(\nu \leqslant k2\), as is immediately evident from the definition of the indices \(k1\) and \(k2\).

Similarly, let \(\omega_{k3}\) and \(\omega_{k4}\) be the two numbers from the sequence with the smallest indices that fall within the interior of the interval \([\omega_a^1, \omega_b^1]\). Let the smaller of these numbers, \(\omega_{k3}\) and \(\omega_{k4}\), be denoted by \(\omega_a^2\), and the larger by \(\omega_b^2\). It is evident that all numbers \(\omega_{\nu}\) from the sequence (\ref{eq13}), for which \(\nu \leqslant k4\), do not fall into the interior of the interval \([\omega_a^2, \omega_b^2]\).

The interval \([\omega_a^\nu, \omega_b^\nu]\) then lies within the interior of all preceding intervals, and it specifically relates to the sequence (\ref{eq13}) in that all numbers \(\omega_{\nu}\), for which \(\nu \leqslant k2\nu\), definitely do not lie within its interior.
  
Since the numbers \(\omega_a^1, \omega_a^2, \omega_a^3, \ldots, \omega_a^\nu, \ldots\) are continually increasing in value while simultaneously being enclosed within the interval \([\alpha, \beta]\), they have a limit that will be denoted by \(A\), such that:
\(A = \lim_{\nu \to \infty} \omega_a^\nu. \)

The same applies to the numbers \(\omega_b^1, \omega_b^2, \omega_b^3, \ldots, \omega_b^\nu, \ldots\), which are continually decreasing and also lie within the interval \([\alpha, \beta]\). Let us denote their limit by \(B\), so that: \( B = \lim_{\nu \to \infty} \omega_b^\nu. \)
The only interesting case to consider is when \(A = B\), to demonstrate that the number:
\(\eta = A = B \)
does not occur in sequence (\ref{eq13}).

\emph{``If it were a member of the sequence, such as the $\nu^{th}$, then one would have: $\eta = \omega_\nu$.
     But the latter equation is not possible for any value of $\nu$ because $\eta$ is in the interior of the interval $[\omega_a^\nu, \omega_b^\nu]$, but $\omega_{\nu}$ lies outside of it."}

The 1879 proof concludes at this point, while the 1873 proof includes an additional explanation:
\emph{``  For all $\nu, \eta\in (\omega_a^\nu, \omega_b^\nu)$ but $\omega_\nu \notin (\omega_a^\nu, \omega_b^\nu)$. Therefore, $\eta$ is a number in $[\alpha,\beta]$ that is not contained in (\ref{eq13})"} \cite{original,originala}.

In the concluding arguments of both proofs, two fundamental issues arise: the technical problem that the proposed construction of a number different from all listed numbers cannot be completed, and the logical problem that assumes if a number does not have an assigned specific \(n\)-th place in the sequence, it is not in the sequence. It will be demonstrated that the number proposed to be different from all others cannot be constructed as intended. What is actually constructed is the number \(\eta\), which is not distinct from other numbers. The construction requires an additional step to differentiate it from the listed numbers, but as it will be shown this final step cannot be completed as proposed. Also, the constructed number \(\eta\) does not have a specific \(n\)-th place in the sequence due to the method of its construction; however, this does not imply that it is not part of the sequence. Let us elaborate on both of these points.

The technical problem lies in the statement: 
\emph{``But the latter equation is not possible for any value of $\nu$ because $\eta$ is in the interior of the interval $[\omega_a^\nu, \omega_b^\nu]$, but $\omega_{\nu}$ lies outside of it.
For all $\nu, \eta\in (\omega_a^\nu, \omega_b^\nu)$ but $\omega_\nu \notin (\omega_a^\nu, \omega_b^\nu)$. Therefore, $\eta$ is a number in $[\alpha,\beta]$ that is not contained in (\ref{eq13})"} \cite{original,originala}.

The problem arises because the constructed number \(\eta = A = B\) is actually from the sequence (\ref{eq13}), which can be easily proven using the NIP. The subsequent step, which aims to distinguish \(\eta\) from the other sequence numbers by requiring that it lies within the open interval \((\omega_a^\nu, \omega_b^\nu) \subset [\omega_a^\nu, \omega_b^\nu]\), 
which eliminates the endpoints \(\omega_a^\nu\) and \(\omega_b^\nu\) in the limit when \(A = B\), cannot be completed.

The sequence elements \(\omega_\nu\) are organized to serve as endpoints in the intervals \([\omega_a^\nu, \omega_b^\nu]\). These intervals form a nested sequence: 
\begin{equation}
[\omega_a^1, \omega_b^1]\supset [\omega_a^2, \omega_b^2]\supset [\omega_a^3, \omega_b^3]\supset \ldots [\omega_a^\nu, \omega_b^\nu]\ldots
\label{eq14}
\end{equation}
that converges to the zero size interval \([A, B]\), which is the number \(\eta = A = B\) defined by that sequence. Each interval of this infinite sequence of nested, closed intervals is nonempty and contains \(\omega_\nu\) numbers from sequence (\ref{eq13}), because that sequence is dense everywhere, as assumed by the construction. According to NIP the intersection of this sequence must contain a single number, which must be from sequence (\ref{eq13}) since no other numbers are present in the sequence (\ref{eq14}) intervals than those from sequence (\ref{eq13}). Consequently, the number \(\eta = A = B\) in the intersection of sequence (\ref{eq14}) is also a member of sequence (\ref{eq13}). 
The sequence (\ref{eq14}) concludes when the endpoints \(\omega_a^\nu\) and \(\omega_b^\nu\) converge, at which point \(A = B\), making \(\eta\) the sole element in all intervals of the sequence. Thus, \(\eta\) exists as an element of both sequences (\ref{eq14}) and (\ref{eq13}), because all numbers in sequence (\ref{eq14}) are derived from sequence (\ref{eq13}) and in each interval \([\omega_a^\nu, \omega_b^\nu]\) by construction there are \(\omega_\nu\) numbers. It is one of the numbers \(\omega_\nu\) from sequence (\ref{eq13}), defined when ultimately \(\omega_a^\nu = \omega_b^\nu\) at the limit \(A=B\).

There is a misstep in proving that all numbers \(\omega_\nu\) are eliminated, leading to the incorrect conclusion that the number \(\eta\) is not included in sequence (\ref{eq13}).

The problem with eliminating all numbers from sequence (\ref{eq13}) using open intervals \((\omega_a^\nu, \omega_b^\nu)\) is as follows: For any finite \(\nu\), no interval \((\omega_a^\nu, \omega_b^\nu)\) can exclude all \(\omega_\nu\) while containing the number \(\eta\), due to the density of sequence (\ref{eq13}). This density implies that any interval \((\omega_a^\nu, \omega_b^\nu)\) contains infinitely many numbers \(\omega_\nu\) from sequence (\ref{eq13}).
Therefore, the sequence numbers \(\omega_\nu\) must be eliminated only in the limit as \(\nu\) approaches infinity. This is proposed to be achieved by designating all numbers \(\omega_\nu\) from sequence (\ref{eq13}) as the endpoints of the nested intervals \([\omega_a^\nu, \omega_b^\nu]\) in sequence (\ref{eq14}), and then eliminating them through the open intervals \((\omega_a^\nu, \omega_b^\nu)\), which are structured to exclude all endpoints \(\omega_a^\nu\) and \(\omega_b^\nu\), and thus all numbers \(\omega_\nu\) in the sequences (\ref{eq13}) and (\ref{eq14}).

The claim is that in the limit as \(\nu \to \infty\), where \(\omega_a^\nu = A\) and \(\omega_b^\nu = B\), and when \(A = B\), there exists an open interval \((A, B)\) that eliminates all endpoints \(\omega_a^\nu\), \(\omega_b^\nu\), and is not empty, containing the number \(\eta=A = B\). To quote:
\emph{``because $\eta$ is in the interior of the interval $[\omega_a^\nu, \omega_b^\nu]$, but $\omega_{\nu}$ lies outside of it.
For all $\nu, \eta\in (\omega_a^\nu, \omega_b^\nu)$ but $\omega_\nu \notin (\omega_a^\nu, \omega_b^\nu)$."}

However, as explained above, the process concludes when the intervals \([\omega_a^\nu, \omega_b^\nu]\) converge to a single point at \([A, B]\), forming a closed, zero-size interval that represents the intersection of sequence (\ref{eq14}), a singleton, a single number. This number must reside within all intervals of sequence (\ref{eq14}) and is derived from sequence (\ref{eq13}).

The final step in the construction, required to make the number \(\eta\) distinct from all listed numbers \(\omega_\nu\), cannot be completed. This step involves using open intervals \((\omega_a^\nu, \omega_b^\nu)\) to eliminate the endpoints of all closed intervals \([\omega_a^\nu, \omega_b^\nu]\) and recognize within that open interval \((\omega_a^\nu, \omega_b^\nu)\) the number \(\eta\). While this procedure can be performed at any finite step \(\nu\), it cannot be executed in the limit as \(\nu \to \infty\). At this limit, it would require that the intersection \([A, B]\) includes an open interval \((A, B)\) that eliminates all \(\omega_\nu\) by removing all endpoints \(\omega_a^\nu\), \(\omega_b^\nu\), and that this open interval also contains the number \(\eta = A = B\), purported to be distinct from all \(\omega_\nu\).

However, according to the NIP, the intersection \([A, B]\), which is a singleton, a single number, cannot encompass an open interval \((A,B)\). Consequently, there can be no open interval \((A,B)\) that includes a number \(\eta=A=B\) different from all numbers in sequence (\ref{eq13}). Therefore, the constructed number \(\eta\), is no different from all sequence's numbers.

The logical problem arises from the statement: \emph{``\(\eta = A = B\) does not occur in sequence (\ref{eq13}). If it were a member of the sequence, such as the \(\nu^{th}\), then one would have: \(\eta = \omega_\nu\)."} 

The problem is the misconception that the inability to assign a specific \(\nu\)-th place in a sequence to a number implies that the number is not in the sequence. This misconception stems from how the number \(\eta\) is defined. It is identified by traversing an infinite number of sequence (\ref{eq14}) elements, which was necessary to eliminate all \(\omega_\nu\) numbers. While \(\eta\) can still be one of the \(\omega_\nu\) numbers, a specific \(\nu\)-th position cannot be assigned due to the method of its definition. To elucidate this point, let us consider the following example, which will demonstrate the incorrectness of this logical conclusion as applied in the proofs.

Consider a sequence of all natural numbers: \(1, 2, 3, \ldots\). Begin by rearranging it, moving the number one to the second position, resulting in \(2, 1, 3, \ldots\). Then rearrange it again, this time placing number one in the \(3^{rd}\) position: \(2, 3, 1, \ldots\). Continue this process of rearrangement an infinite number of times, which is feasible since the sequence is infinite. After an infinite number of rearrangements, it becomes impossible to specify the exact \(\nu\)-th position of the number one. However, it remains within the sequence; it has never been removed, only repositioned.

The situation with the missing number \(\eta\), which is defined using an infinite procedure, is similar. After traversing all elements of the sequence (\ref{eq14}), infinitely many of them, necessary to eventually eliminate all endpoints of the closed intervals \([\omega_a^\nu, \omega_b^\nu]\), which can only occur in the limit as \(\nu \rightarrow \infty\) when \(\omega_a^\nu = \omega_b^\nu\), \(\eta\) is finally identified. Because \(\eta\) emerges through an infinite process, going over an infinite number of  sequence elements \(\omega_\nu\), its \(\nu\)-th position cannot be specified. However, this does not imply that it is not a member of the sequence, as demonstrated by the given example of the permuted sequence of natural numbers.

One or sometimes both of these problems occur in all other proofs of the non-denumerability of real numbers. For instance, in the diagonal argument, the issue arises because the proposed construction cannot be completed in the limit, which is where the diagonal number is defined. Consequently, it becomes impossible to construct a number that is different from all listed sequence numbers, a point that will be demonstrated.

{\bf The diagonal argument}, constructs the number \(\eta\), designed to differ in the \(n\)-th decimal place from the \(n\)-th listed number \(\omega_n\), implying that the list of numbers \(\omega_n\) is incomplete. This concept can also be interpreted through nested intervals. To facilitate this, we redefine the diagonal number \(\eta\) and the sequence of listed numbers \(\omega_n\) using a framework of nested intervals instead of the traditional decimal expansion. In this representation, the diagonal argument ensures that at each finite step, \(\eta\) differs from each \(\omega_n\) by having sequences representing \(\eta\) and \(\omega_n\) differ in the \(n\)-th subinterval.

As will be demonstrated, the diagonal argument, when represented through nested intervals, conflicts with the common interpretation that the number \(\eta\), which diverges from each \(\omega_n\) at finite stages, can also be constructed in the limit, which must be confronted when considering all listed sequence numbers \(\omega_n\). The number \(\eta\) is fully defined only in the limit, after all its decimal places in the decimal expansion or all intervals in the nested interval representation are determined, after traversing through all listed numbers.
At any finite stage, the construction is feasible since each interval \(I_n\) that represents \(n-th\) decimal place contains two disjoint subintervals, one \(H_n\) assigned to \(\eta\) and the other \(\Omega_n\) to \(\omega_n\). However, in the limit, the proposed construction of the diagonal number \(\eta\) cannot be completed because it would contradict the NIP. This is because, in the limit, the intersection of intervals \(I_n\) reduces to a singleton, a single number. Consequently, in the limit, it becomes impossible for there to be two distinct subintervals, one \(H_n\) representing \(\eta\) and another \(\Omega_n\) representing \(\omega_n\), to ensure that \(\eta\) differs from all \(\omega_n\). Instead, constructed is number \(\eta\) that is not different from all other listed numbers. The final step, needed to make it different from all other numbers, which by the construction requires assigning \(\eta\) a subinterval \(H_n\) that is distinct from the subinterval \(\Omega_n\) assigned to number \(\omega_n\), cannot be completed.
Thus, the proposed construction of \(\eta\), achievable at all finite stages, cannot be completed because both \(H_n\) and \(\Omega_n\) represent the same decimal and therefore belong to the same intervals \(I_n\), which at the limit have an intersection that is a singleton, not allowing disjunctive intervals \(H_n\) and \(\Omega_n\) as required by the construction. 

The diagonal argument assumes that the conceptual process of constructing the number \(\eta\) can theoretically be extended indefinitely. However, when this concept is translated into a practical framework, such as interval representation, significant implications arise at the limit, leading to the conclusion that such construction cannot be completed. The consideration of the limit is essential because the number \(\eta\) is fully defined only in this context.
While the diagonal argument, in its original formulation, ensures divergence at each finite stage without explicit concern for the limit, adopting an interval-based approach underscores the need to reconcile this finite-stage divergence with what occurs in the limit.
Let us explore this explicitly.

Numbers \(\eta\) and \(\omega_n\) in decimal representation can be written as 
\begin{equation}
   0.a_1 a_2 a_3 \cdots a_n \cdots =\sum_{n=1}^{\infty}a_n10^{-n} 
\label{eq15}
\end{equation}
where $a_n$ are 0,1,2,3,4,5,6,7,8,9. 

In the interval representation, the decimal places \(a_n\) are represented by nested intervals \(I_n = [\alpha_n, \beta_n]\), where the sizes of the intervals \(I_n = |\alpha_n - \beta_n| = 10^{-n}\) correspond to the \(n\)-th decimal place positions. These are specified by decimals \(a_n\) in (\ref{eq15}), ensuring that \(a_n \in [\alpha_n, \beta_n]\). Each interval \([\alpha_n, \beta_n]\) represents the \(n^{th}\) decimal in (\ref{eq15}), with \([\alpha_n, \beta_n] \supset [\alpha_{n+1}, \beta_{n+1}]\).

The numbers \(\eta\) and \(\omega_n\) are then defined through an infinite sequence of nested intervals
\begin{equation}
          I_1 \supset  I_2 \supset  I_3 \supset. . . I_n \supset I_{n+1} \supset  . . .     
\label{eq16}
\end{equation}
i.e., they are intersections of these \(I_n\) intervals
\begin{equation}
\bigcap_{n \in \mathbb{N}} I_n. 
\label{eq17}
\end{equation}

The number \(\eta\) constructed by the diagonal argument is defined such that each decimal \(a_n\) differs from the \(n^{th}\) decimal place of the corresponding number \(\omega_n\) in the list. Thus, \(a_n \neq \omega_{nn}\), where \(\omega_{nn}\) represents the \(n^{th}\) digit of the number \(\omega_n\) that occupies the \(n^{th}\) position in the listed sequence of \(\omega_n\) numbers.

In the nested interval representation, the equivalent expression stipulates that at each finite step \(n\), the process involves considering the intersection (\ref{eq17}) of the first \(n\) intervals \(I_n\) representing the first \(n\) decimals of the numbers \(\eta\), and \(\omega_n\). Within these intersections, distinct subintervals \(H_{n}\) for \(\eta\) and \(\Omega_{n}\) for \(\omega_n\) are selected to ensure divergence. Thus the divergence is created by choosing disjunctive subintervals to guarantee that \(\eta\) and \(\omega_n\) differ at the \(n\)-th stage, effectively making the numbers different at each \(n\)-th decimal place or \(n\)-th interval, because they will have different intersections (\ref{eq17}):
\begin{equation}
\eta=\bigcap_{i \in \mathbb{N}} H_i, \quad\quad \mathrm{and}
%\eta=\bigcap_{i \in \mathbb{N}} H_i \quad \mathrm{and} \quad \omega_n=\bigcap_{i \in \mathbb{N}} \Omega_i. 
\label{eq18}
\end{equation}
\begin{equation}
\omega_n=\bigcap_{i \in \mathbb{N}} \Omega_i. 
%\eta=\bigcap_{i \in \mathbb{N}} H_i \quad \mathrm{and} \quad \omega_n=\bigcap_{i \in \mathbb{N}} \Omega_i. 
\label{eq19}
\end{equation}
The diagonal argument assumes that the process of creating disjunctive subintervals \(H_i\) and \(\Omega_i\) can continue indefinitely, theoretically allowing for the construction of the number \(\eta\) different from all listed numbers \(\omega_n\). However, this seemingly straightforward assumption contradicts the NIP. At each finite stage in sequence (\ref{eq16}), every interval \(I_n = [\alpha_n, \beta_n]\) has a finite size, with \(\alpha_n < \beta_n\) consistently satisfied. This ensures that each interval \(I_n\) contains at least two disjunctive subintervals \(H_n\) and \(\Omega_n\), necessary to construct number $\eta$ such that it will be different from listed numbers $\omega_n$. However, this intuitively natural requirement cannot be met when applied to the entire sequence (\ref{eq16}) as its intersection interval (\ref{eq17}), in the limit, has zero size.

     In the limit, as \(n \rightarrow \infty\), \(\alpha_n \rightarrow A\) and \(\beta_n \rightarrow B\), and \(A = B\), the intersection (\ref{eq17}) of all intervals \(I_n = [\alpha_n, \beta_n]\), i.e., the intersection interval \([A, B]\), has zero size and 
contains only a single number.   So, there is no possibility that in the intersection interval will exist two disjunctive subintervals \(H_n\) and \(\Omega_n\), which is required by diagonal argument to make the intersection (\ref{eq18}) representing number $\eta$ different from the intersections (\ref{eq19}) representing a listed number $\omega_n$. 

   The intersection theorem dictates that sequence (\ref{eq16}) has a zero-size intersection with only one number in it. Intervals \(H_n\) representing decimals \(a_n\) in number \(\eta\) and intervals \(\Omega_n\) representing \(\omega_{nn}\) in \(\omega_n\) for every \(n\) belong to the same interval \(I_n=[\alpha_n, \beta_n]\) because they represent the same \(n\)-th decimal position. Therefore, they must lie within the same intersection (\ref{eq17}) of intervals \(I_n\) at each step for any \(n\).
   For every \(n\) in the construction of numbers \(\eta\) and \(\omega_n\), the intersection of the first \(n\) intervals \([\alpha_n, \beta_n]\) has non-zero size and contains disjunctive subintervals \(H_n\) and \(\Omega_n\), allowing \(\eta\) to differ from the first \(n\) listed numbers \(\omega_n\) in their \(n\)-th decimal places. However, at these finite stages, number \(\eta\) is not yet defined or constructed, regardless of how large \(n\) is.
    Number \(\eta\) is only defined or constructed in the limit as \(n \rightarrow \infty\), when \(\alpha_n \rightarrow A\) and \(\beta_n \rightarrow B\). But in this limit, the intersection \([A, B]\) is a singleton containing only one number. The requirement at each \(n\)-th step to select two disjunctive subintervals \(H_n\neq \Omega_n\) that make sequence intersections (\ref{eq18}) and (\ref{eq19}) different, cannot be fulfilled in the limit because the intersection (\ref{eq17}) becomes a singleton, a single number. So, the number $\eta$, that is different from all listed numbers cannot be constructed since in the limit do not exist two disjunctive subintervals, required to make the number $\eta$ different from all listed $\omega_n$ numbers. 

\section{Conclusion}

Discussed are differences between various models for generating sets of numbers and how the inclusion of specific axioms influences the cardinality of these generated sets. Considered are some open  questions and discrepancies among the current leading models, particularly those represented by V - ultimate L and models based on Martin's axioms. Additionally, two explicit constructions of generated sets of real numbers are analyzed in detail.

It is demonstrated that expanding the \(ZF\)  framework with the Axiom of Choice exposes inherent uncertainties in determining the cardinality of infinite sets. The constructed sets,  which are both everywhere dense and Dedekind complete yet countable, align with the  Löwenheim-Skolem theorem. The role of \(AC\) is emphasized in facilitating the denumerability of these generated complete sets, while tensions are identified between the Nested Interval Property and classical proofs of non-denumerability.

The second generated set serves as both a sequence representing all real numbers and the first explicit example of Wang's constructive model. This dual role highlights the interplay between classical and constructive approaches to generating the continuum, offering a bridge between theoretical frameworks and concrete constructions.

These findings underscore the intricate relationship between axiomatic extensions, foundational properties such as completeness and density, and the nature of infinite sets. They provide a foundation for further exploration of the implications of specific axioms on set cardinality and the continuum hypothesis in modern set theory.

\vskip 0.2cm

{\bf Acknowledgment} 

This work is supported by the National Science Foundation Center of Research Excellence in Science and Technology NSF CREST 1647022, and NSF Research Infrastructure for Science and Engineering NSF RISE 1829245 awards.

\end{document}